\documentclass{amsart}

\usepackage{amssymb,amsfonts, amsmath, amsthm}
\usepackage[all,arc]{xy}
\usepackage{enumerate}
\usepackage{mathrsfs}
\usepackage{mathtools}
\usepackage[mathscr]{euscript}
\usepackage{array}
\usepackage{appendix}

\usepackage{tikz}
\usepackage{tikz-cd}
\usepackage{mathpazo}
\usepackage{textcomp}
\usepackage{bbold}
\usepackage[bbgreekl]{mathbbol}

\usepackage[pagebackref, colorlinks=true, linkcolor=blue, citecolor = red]{hyperref}
 \renewcommand*{\backref}[1]{}
\renewcommand*{\backrefalt}[4]{({%
		\ifcase #1 Not cited.%
		\or On p.~#2%
		\else On pp.~#2%
		\fi%
	})}

\usepackage[nameinlink,capitalise,noabbrev]{cleveref}
\crefname{subsection}{Subsection}{Subsection}

\usetikzlibrary{patterns}

\usepackage{geometry}
\usepackage{todonotes}

\newcommand{\s}{\mathscr{S}}
\renewcommand{\ss}{s\mathscr{S}}
\newcommand{\sss}{ss\mathscr{S}}

\newcommand{\C}{\mathscr{C}}
\newcommand{\D}{\mathscr{D}}
\newcommand{\E}{\mathscr{E}}

\newcommand{\M}{\mathcal{M}}

\newcommand{\R}{\mathscr{R}}

\newcommand{\sSet}{\mathrm{s}\mathscr{S}\mathrm{et}}

\newcommand{\kX}{\mathfrak{X}}
\newcommand{\kY}{\mathfrak{Y}}
\newcommand{\kV}{\mathfrak{V}}
\newcommand{\kW}{\mathfrak{W}}
\newcommand{\kD}{\mathfrak{D}}

\newcommand{\kL}{\mathfrak{L}}
\newcommand{\kF}{\mathfrak{F}}
\newcommand{\kZ}{\mathfrak{Z}}
\newcommand{\kC}{\mathfrak{C}}
\newcommand{\ks}{\mathfrak{s}}
\newcommand{\kt}{\mathfrak{t}}
\newcommand{\ki}{\mathfrak{i}}
\newcommand{\kj}{\mathfrak{j}}
\newcommand{\kf}{\mathfrak{f}}
\newcommand{\kx}{\mathfrak{x}}
\newcommand{\ky}{\mathfrak{y}}
\newcommand{\kl}{\mathfrak{l}}
\newcommand{\kc}{\mathfrak{c}}
\renewcommand{\kt}{\mathfrak{t}}
\newcommand{\comp}{\mathrm{comp}}

\newcommand{\id}{\mathrm{id}}
\newcommand{\set}{\mathscr{S}\text{et}}
\newcommand{\cat}{\mathscr{C}\text{at}}

\newcommand{\Cone}{\mathscr{C}\mathrm{one}}

\newcommand{\Sub}{\mathscr{S}\mathrm{ub}}

\newcommand{\Hom}{\mathrm{Hom}}
\newcommand{\Map}{\mathrm{Map}}

\newcommand{\Fun}{\mathrm{Fun}}
\newcommand{\map}{\mathrm{map}}
\newcommand{\comma}{,}

\newcommand{\LFib}{\mathscr{L}\mathscr{F}\mathrm{ib}}
\newcommand{\RFib}{\mathscr{R}\mathscr{F}\mathrm{ib}}
\newcommand{\newLFib}{\mathfrak{L}\mathfrak{F}\mathrm{ib}}
\newcommand{\coCartFib}{\mathrm{co}\mathfrak{C}\mathrm{art}\mathfrak{F}\mathrm{ib}}
\newcommand{\CSS}{\mathscr{C}\mathscr{S}\mathscr{S}}

\newcommand{\Diag}{\mathscr{D}\mathrm{iag}}
\newcommand{\Inc}{\mathscr{I}\mathrm{nc}}
\newcommand{\Und}{\mathscr{U}\mathrm{nd}}
\renewcommand{\Row}{\mathscr{R}\mathrm{ow}}
\newcommand{\infref}[1]{\cref{subsec:infty categories}(\ref{#1})}
\newcommand{\covref}[1]{\cref{subsec:covar}(\ref{#1})}
\newcommand{\reecovref}[1]{\cref{subsec:reedy covar}(\ref{#1})}
\newcommand{\whyref}[1]{\cref{subsec:why}(\ref{#1})}

\DeclareMathAlphabet{\mathbbe}{U}{bbold}{m}{n}

\def\DDelta{{\mathbbe{\Delta}}}
\renewcommand{\1}{\mathbbe{1}}
\renewcommand{\2}{\mathbbe{2}}
\newcommand{\3}{\mathbbe{3}}

\newcommand{\adjun}[4]{
\begin{tikzcd}[row sep=0.5in, column sep=0.5in]
 #1  \arrow[r, shift left=1.8, "#3", "\bot"'] \pgfmatrixnextcell
 #2 \arrow[l, shift left=1.8, "#4"] 
\end{tikzcd}
}

\newcommand{\simpset}[7]{

 \begin{tikzcd}[row sep=0.5in, column sep=0.5in]
   #1 \arrow[r, shorten >=1ex,shorten <=1ex]
   \pgfmatrixnextcell #2 
   \arrow[l, shift left=1.2, "#5"] \arrow[l, shift right=1.2, "#4"'] 
   \arrow[r, shift right, shorten >=1ex,shorten <=1ex ] \arrow[r, shift left, shorten >=1ex,shorten <=1ex] 
   \pgfmatrixnextcell #3 
   \arrow[l] \arrow[l, shift left=2, "#7"] \arrow[l, shift right=2, "#6 "'] 
   \arrow[r, shorten >=1ex,shorten <=1ex] \arrow[r, shift left=2, shorten >=1ex,shorten <=1ex] \arrow[r, shift right=2, shorten >=1ex,shorten <=1ex]
   \pgfmatrixnextcell \cdots 
   \arrow[l, shift right=1] \arrow[l, shift left=1] \arrow[l, shift right=3] \arrow[l, shift left=3] 
 \end{tikzcd}
}

\newtheorem{theone}[equation]{Theorem}
\newtheorem{lemone}[equation]{Lemma}
\newtheorem{propone}[equation]{Proposition}
\newtheorem{corone}[equation]{Corollary}

\theoremstyle{definition}
\newtheorem{defone}[equation]{Definition}

\theoremstyle{remark}
\newtheorem{remone}[equation]{Remark}
\newtheorem{notone}[equation]{Notation}

\newtheoremstyle{TheoremNum}
{}{}              
{\itshape}                      
{}                              
{\bfseries}                     
{.}                             
{ }                             
{\thmname{#1}\thmnote{ \bfseries #3}}
\theoremstyle{TheoremNum}

\numberwithin{equation}{section}

\makeatletter
\def\@seccntformat#1{%
	\expandafter\ifx\csname c@#1\endcsname\c@section\else
	\csname the#1\endcsname\quad
	\fi}
\makeatother

\title{Cartesian Fibrations and Representability}

\author{Nima Rasekh}
\address{{\'E}cole Polytechnique F{\'e}d{\'e}rale de Lausanne, SV BMI UPHESS, Station 8, CH-1015 Lausanne, Switzerland}
\email{nima.rasekh@epfl.ch}
\date{February 2021}

\begin{document}

\begin{abstract} 
 We use the complete Segal approach to the theory of Cartesian fibrations (developed in \cite{rasekh2021cartfibcss}) to define and study representable Cartesian fibrations, generalizing representable right fibrations which have played a key role in $\infty$-category theory. 
 In particular, we give a construction of representable Cartesian fibrations using over-categories and prove the Yoneda lemma for representable Cartesian fibration, which generalizes the established Yoneda lemma for right fibrations. 
 
 We then use the theory of Cartesian fibrations to study complete Segal objects internal to an $\infty$-category. Concretely, we prove the {\it fundamental theorem of complete Segal objects}, which characterizes equivalences of complete Segal objects. Finally we give two application of the results. First, we present a method to construct Segal objects and second we study the representability of the universal coCartesian fibration.
\end{abstract}

\maketitle
\addtocontents{toc}{\protect\setcounter{tocdepth}{1}}

\section{Introduction}\label{Sec Introduction}

\subsection{From Moduli Objects to Representability} \label{subsec:moduli rep}
 The concept of {\it representability} plays an important result in mathematics and goes back at least to Riemann, who is claimed to have coined the term {\it moduli space} with the goal of studying {\it Riemann surfaces} \cite{riemann1857abel,riemann1953gesammelt}. Since then mathematicians have defined objects that classify a certain property and have then studied the object. Further examples in differential geometry includes the {\it classifying space for principal bundles} \cite{milnor1956universalbundlei,milnor1956universalbundleii} and the {\it moduli space of framed manifolds} \cite{randalwilliams2014moduli}. 
 Beside differential geometry moduli spaces are also very prominent in {\it algebraic geometry}: examples include
 {\it Chow varieties} \cite{chevalley1958chowverietyi,chevalley1958chowvarietyii}, {\it Hilbert schemes} \cite{grothendieck1995hilbertschemes}, {\it moduli stack of elliptic curves} \cite{mumford1963ellipticcurves}. The last example also points to the importance of moduli objects in algebraic topology, where beside the moduli stack of elliptic curves \cite{landweberravenelstong1995ellipticcohomology}, the {\it moduli stack of formal groups} \cite{quillen1969fglmu} is used to study homotopy groups of spheres.
 For further discussion of moduli problems in geometry see \cite{benzvi2008moduli}.
 
 The study of moduli objects can be formalized via {\it category theory}. Starting with a category of interest $\C$ we can often describe the desired property as a functor to the category of sets $F: \C^{op} \to \set$. A moduli object in $\C$ is then precisely an object that represents the functor $F$. Hence, the study of moduli objects translates to the representability of functors. 
 In fact we can use this perspective with category theory itself by applying the Yoneda lemma \cite{maclane1998yoneda}. Many categorical notions, such as the existence of limits, colimits and adjunctions can be articulated via universal properties, which in turn can be expressed as a certain set-valued functor. Hence the existence of the desired universal object (the limits or adjunction) translates to the representability of a certain functor. 
 This has been used effectively to study arbitrary categories via the category of sets \cite{maclane1998categories}. 
 We can use this line of thinking to introduce new moduli objects inside categories. One important example is the {\it subobject classifier}, which is the key feature of an {\it elementary topos} \cite{tierney1973elementarytopos}.
 
 We can generalize questions about representability from categories to {\it $\mathscr{V}$-enriched categories}. 
 Now for a given $\mathscr{V}$-enriched category $\C$ we try to address the representability of functors $F: \C^{op} \to \mathscr{V}$.
 This in particular is can be applied to study {\it weighted limits} in enriched categories \cite{kelly1982enriched,riehl2017context}.
 The next logical step would be to generalize questions of representability from enriched categories to {\it $\infty$-categories}.
 Unfortunately, this faces significant challenges, as the non-strict composition in $\infty$-categories prevents us from easily defining representable functors. The solution is to use {\it fibrations}: We describe moduli problems in higher categories by constructing appropriate {\it right fibrations} and then we try to solve them by relating them to {\it representable right fibrations}. Again, we can construct appropriate right fibrations whose solutions give us $\infty$-categorical limits, colimits and adjunctions \cite{lurie2009htt,riehlverity2018elements}. 

\subsection{From Representability to Complete Segal Objects} \label{subsec:rep cso}
 Sometimes the functor we are trying to represent does not naturally take value in sets. A prominent example are {\it cohomology theories}, which are classically defined as functors from the homotopy category of spaces to the category of abelian groups \cite{eilenbergsteenrod1945axioms}. In this case it does not suffice to determine an object to make this functor representable. Rather we need an {\it abelian group object}: a space $A$ along with a morphism $\mu: A \times A \to A$ that satisfies the axioms of an abelian group. 
 We can in fact generalize this argument to any category with finite products. 
 
 A more complicated example is the representability of functors $F:\C^{op} \to \cat$. As categories consist of objects and morphisms we need more data to represent such a functor: Assuming $\C$ is finite complete, we need an object $\mathcal{O}$ (representing the objects), an object $\mathcal{M}$ (representing the morphisms) and various morphisms, such as a composition map $\mu: \mathcal{M} \times_\mathcal{O} \mathcal{M} \to \mathcal{M}$. The data will precisely assemble into an {\it internal category} \cite{roberts2012internalcat}, which also have been studied independent of any representability concerns in a topos theoretic context \cite{johnstone2002elephanti,johnstone2002elephantsii,maclanemoerdijk1994topos}. For example internal categories in the category of sets give us small categories, where as internal categories in the category of small categories give us small {\it double categories} \cite{ehresmann1963doublecat}. The example of internal categories in the category of sets is instructive: by characterizing a small category via a set of objects, a set of morphisms and various diagrams, we could generalize the definition to every finitely complete category. 
 
 We want to generalize this analysis to $\infty$-categories, meaning we want to study the representability of functors $\C^{op}\to \cat_\infty$, where $\C$ is an $\infty$-category with finite limits. Again, it is immediate that one object does not suffice. However, unlike the classical case it also would not suffice to choose two objects, as in an $\infty$-category composition is part of the necessary data. Fortunately, a solution has already been worked out for us: {\it complete Segal spaces} \cite{rezk2001css}. A complete Segal space is a simplicial object in spaces, which satisfies two conditions: the {\it Segal condition} and the {\it completeness conditions}, both of which can be characterized via finite limit diagrams. Moreover a complete Segal space is a model for an $\infty$-category \cite{bergner2010survey,bergner2018book,joyaltierney2007qcatvssegal,toen2005unicity}.
 Using complete Segal spaces as motivation, we can hence define a {\it complete Segal object} inside a finitely complete $\infty$-category as a simplicial object that satisfies the Segal and completeness condition and choose that as our model of internal $\infty$-category. \footnote{Complete Segal objects have already been studied under the name {\it Rezk object} in \cite{riehlverity2017inftycosmos} (with other applications).}
 
 Having understood the role of complete Segal objects as internal $\infty$-categories, we can now study the representability of functors valued in $\infty$-categories. Again, we need to translate our discussion from functors to fibrations: this time we get {\it Cartesian fibrations} \cite{lurie2009htt}. 
 Hence, we want to understand the relationship between {\it representable Cartesian fibrations} and complete Segal objects. 
 This is the goal of this work.
 
 \subsection{Simplicial Objects and Representable Reedy Left Fibrations}
 In fact we will study this relationship as a special case of results that are more general in two respects.
 First, rather than using any specific model of $\infty$-categories (such as {\it quasi-categories} \cite{joyal2008notes,joyal2008theory}), we take a completely model independent approach: {\it a theory of $\infty$-categories with fibration structure} (\cref{sec:infty}). 
 This includes a model structure $\M$, which comes with an underlying {\it $\infty$-cosmos} \cite{riehlverity2017inftycosmos,riehlverity2018elements} (\cref{subsec:infty categories}), along with a model structure of left (right) fibrations (\cref{subsec:covar}). This will in particular include quasi-categories and complete Segal spaces as examples. (\cref{subsec:examples}). 
  
 Next, we will in fact start by studying the representability of fibrations that correspond to functors valued in simplicial spaces, so-called {\it Reedy left (right) fibrations} \cite{rasekh2021cartfibcss, rasekh2021cartfibmarkedvscso}, which exists in every $\infty$-category with fibration structure (\cref{subsec:reedy covar}). We will prove that for every cosimplicial object $\kx^\bullet: \DDelta \to X$, there exists a Reedy left fibration, we denote $X_{\kx^\bullet/} \to X$ that we call the {\it representable Reedy left fibration} (\cref{def:rep reedy left fib}). Moreover, if $\C$ is an $\infty$-category and $\kx^\bullet: \DDelta \to \C$, then we can give an explicit construction of the representable Reedy left fibration $\C_{\kx^\bullet/} \to \C$ via {\it cocones} (\cref{prop:rep left fib higher cat}), which in particular means that the fiber over an object $\ky$ in $\C$ is the simplicial space
 \begin{center}
 	\simpset{\map_\C(\kx^0, \ky)}{\map_\C(\kx^1, \ky)}{\map_\C(\kx^2, \ky)}{}{}{}{}.
 \end{center}
 Moreover, these representable fibrations satisfy the {\it Yoneda lemma} (\cref{the:yoneda Lemma for cosimplicial objects}): For two cosimplicial objects $\kx^\bullet, \ky^\bullet: \DDelta \to\C$ we have an equivalence 
 $$\Map_{/\C}(\C_{\kx/},\C_{\ky/}) \xrightarrow{ \ \simeq \ } \map_{cos\C}(\ky^\bullet,\kx^\bullet).$$
 Having these results we can move on to representable Cartesian fibrations and complete Segal objects. We prove that a given simplicial object $\kx^\bullet: \DDelta^{op} \to \C$, $\kx^\bullet$ is a complete Segal object if and only if $\C_{/\kx_\bullet}$ is a Cartesian fibration (\cref{prop:segal obj vs Segal Cart}). We use this perspective to study the category theory of complete Segal objects, in particular proving the {\it fundamental theorem of complete Segal objects}, which characterizes their weak equivalences (\cref{the:fundamental theorem cso}).
 
\subsection{Why Representable Cartesian Fibrations and Complete Segal Objects?} \label{subsec:why}
 What are the benefits of studying representable Cartesian fibrations, complete Segal objects and their interaction?
 \begin{enumerate}
 	\item \label{item:constructing cso} {\bf Constructing Internal $\infty$-Categories:} Let $\C$ be a finite complete category and $\mathcal{O}$ and $\mathcal{M}$ be two objects in $\C$ such that for every object $c$ in the $\C$, the sets $\Hom(c,\mathcal{O}),\Hom(c,\mathcal{M})$ are sets of objects and morphisms of a category, natural in $c$. Then $\mathcal{O},\mathcal{M}$ are in fact part of an internal category in $\C$. Indeed, we can use the Yoneda lemma to translate the external maps (such as source, target and composition) to the objects $\mathcal{O},\mathcal{M}$. 
 	
 	We would like a similar result in the $\infty$-categorical world, meaning whenever we have two objects $\kW_0,\kW_1$ in a finitely complete $\infty$-category $\C$ and the structure of an $\infty$-category (here complete Segal space) on the spaces $\map_\C(\kx,\kW_0)$,$\map_\C(\kx,\kW_1)$ then it in fact gives us a complete Segal object in $\C$. Unfortunately, there is no direct way to prove such result. We can use pullbacks to construct each level $\kW_n \simeq \kW_1 \times_{\kW_0} ... \times_{\kW_0} \kW_1$, however, we cannot simply construct the maps between the $\kW_n$ by hand as there are an infinite number of them and they need to be coherent.
 	
 	The solution is to translate our problem into the language of fibrations and solve it there using a powerful strictification result, which we do in \cref{subsec:strictify}.
 	
  	\item \label{item:universal fib} {\bf Representing the Universal coCartesian Fibration:} Applying the logic of representability (from \cref{subsec:moduli rep}) on left fibrations of $\infty$-categories, we would like to know there exists a {\it universal left fibration} that classifies all other left fibrations (via pullback). Fortunately, this question has already been addressed by many authors: the universal left fibration is the projection map $\s_* \to \s$ from pointed spaces to spaces, which is in fact a representable left fibration, represented by the point \cite[Lemma 2.2.4]{kazhdanvarshvsky2014yoneda} \cite[Corollary 5.2.8]{cisinski2019highercategories}. 
  	
  	We could ask a similar question regarding coCartesian fibrations. Does there exists a coCartesian fibration that classifies all other coCartesian fibrations? We can give an abstract answer \cite[Subsection 3.3.2]{lurie2009htt} , however, we like to be able to explicitly describe the resulting Cartesian fibration, similar to how we described the universal left fibration. We cannot describe it using a representing object, however, in \cref{subsec:universal} we will explain how the {\it universal coCartesian fibration} is represented by a complete Segal object. 
  	Note this has also been discussed extensively in \cite[Example 3.26]{stenzel2020comprehension}.
  	
  	\item \label{item:target fib} {\bf Representing the Target Fibration and Topos Theory:}
  	An elementary topos $\E$ is characterized via a {\it subobject classifier}, which represents the functor $\Sub:\E^{op} \to \set$ that takes an object $c$ in $\E$ to the full subcategory of monomorphisms in  $\E_{/c}$ \cite{maclanemoerdijk1994topos}. 
  	We want to generalize this definition to the $\infty$-categorical setting, which necessitates representing an appropriate $\infty$-categorical generalization of $\Sub$. For a finitely complete $\infty$-category $\E$, the fibration that corresponds to the functor with value $\E_{/-}$
  	is given by the {\it target Cartesian fibration} from the arrow category $t:\E^\2 \to \E$. 
  	
  	As $t: \E^\2\to \E$ is a Cartesian fibration its representability has not been studied, rather the focus has been on the sub-fibration of $\mathscr{O}_\E\to \E$, that corresponds to the functor with value the maximal subgroupoid of $\E_{/-}$. 
  	This is in fact a right fibration, and its representability has been studied extensively, with representing objects known as  
  	{\it object classifiers} or {\it universes}, and with many applications in {\it $\infty$-topos theory} \cite{lurie2009htt} and {\it homotopy type theory} \cite{shulman2019inftytoposunivalent}. 
  	
  	Developing a functioning theory of representable Cartesian fibrations allows us to study the representability of the target fibration itself, rather than any sub-fibration, which results in generalizations of universes which are expected to play important roles in higher topos theory \cite{rasekh2018elementarytopos}. 
  	Hence studying the representability of important Cartesian fibrations can enable us define new concepts and universal properties in $\infty$-category theory, with the example of target fibrations hopefully being only the first of many.
\end{enumerate}

\subsection{Relation to Other Work} 
This paper is the third part of a three-paper series which introduces the bisimplicial approach to Cartesian fibrations: 
\begin{enumerate}
	\item {\bf Cartesian Fibrations of Complete Segal Spaces} \cite{rasekh2021cartfibcss}
	\item {\bf Quasi-Categories vs. Segal Spaces: Cartesian Edition} \cite{rasekh2021cartfibmarkedvscso}
	\item {\bf Cartesian Fibrations and Representability} 
\end{enumerate}
In particular, the first paper introduced the complete Segal approach to Cartesian fibrations in the complete Segal space setting. The second paper translated it to the setting of quasi-categories and proved it is equivalent to marked simplicial set approach. Those papers serve as a foundational backbone of this third paper which applies it to representable Cartesian fibrations. 

It should be noted that after this work first appeared Stenzel \cite[Subsection 4.1]{stenzel2020comprehension} studied a lot of the same topics, largely independently, using Cartesian fibrations in the sense of Lurie \cite{lurie2009htt}.

\subsection{Acknowledgments} \label{Subsec Acknowledgements}
I want to thank my advisor Charles Rezk for suggesting the study of internal $\infty$-categories via complete Segal objects.
I also want to thank Raffael Stenzel for studying representable Cartesian fibration from the perspective of the comprehension construction in \cite{stenzel2020comprehension} and an excellent talk in the Masaryk University Algebra Seminar, which has motivated \cref{subsec:universal}.
 
\section{A Theory of \texorpdfstring{$\infty$}{oo}-Categories with Fibration Structure} \label{sec:infty}
In this section we introduce an axiomatic structure we call a {\it theory of $\infty$-categories with fibration structure} that will serve as our $(\infty,1)$-categorical foundations in the coming sections (with the sole exception of \cref{rem:segal also works} where we mention Segal spaces).

A theory of $\infty$-categories with fibration structure consists of following data:
\begin{enumerate}
	\item[(I)] A left proper combinatorial model category compatible with Cartesian closure $\M$ \cite{hovey1999modelcategories,hirschhorn2003modelcategories}. 
	\item[(II)] Two Quillen equivalences 
	\begin{equation} \label{eq:adjunctions}
		\begin{tikzcd}[column sep=0.5in]
		 \ss^{CSS} \arrow[r, shift left=1.8, "\Diag", "\bot"'] & \M\arrow[l, shift left=1.8, " \Und"] \arrow[r, shift left=1.8, "\Inc", "\bot"'] & \ss^{CSS} \arrow[l, shift left=1.8, "\Row"]
		\end{tikzcd} 
	\end{equation}
 	where the left adjoint $\Diag$ preserves binary products of $F(n)$ (\infref{item:fn}).
	Here the category of simplicial spaces has the model structure for complete Segal spaces \cite[Theorem 7.2]{rezk2001css}.
	\item[(III)] For each object $X$ in $\mathcal{M}$ there exists two simplicial combinatorial left proper left Bousfield localization model structures on the over-category $\M_{/X}$, called the {\it covariant} and {\it contravariant} model structures, such that the induced adjunction 
	\begin{equation} \label{eq:relative adjunctions}
		\adjun{(\M_{/X})^{cov}}{(\ss_{/\Inc(X)})^{cov}}{\Inc}{\Row}, \
		\adjun{(\M_{/X})^{contra}}{(\ss_{/\Inc(X)})^{contra}}{\Inc}{\Row}
	\end{equation}
	are Quillen equivalences. Here the left hand category of simplicial spaces the {\it covariant model structure} and the right hand category of simplicial spaces has the {\it contravariant model structure}, both of which have been defined in \cite{rasekh2017left}. 
\end{enumerate}

For the remainder of this section we establish some important notation and results regarding $\M$.

\subsection{Model Structure for \texorpdfstring{$\infty$}{oo}-Categories} \label{subsec:infty categories}
We want to establish some results regarding the theory of $\infty$-categories.
\begin{enumerate}
	\item \label{item:gen} We call the objects in $\M$ the {\it generic objects} and denote them with letters $X, Y, Z$. 
	\item The existence of the Quillen equivalence $(\Diag,\Und)$ is precisely the condition of being a {\it theory of $(1,\infty)$-categories} in the sense of \cite{toen2005unicity}. 
	\item The full subcategory of fibrant objects $\M^f$ is an $\infty$-cosmos as defined in \cite[Definition 2.1.1]{riehlverity2017inftycosmos}. 
	This follows from the fact that the adjunction 
	\begin{center}
		\begin{tikzcd}[column sep=0.7in]
			\sSet \arrow[r, shift left=1.8, "\Diag\circ p_1^*", "\bot"'] & \M  \arrow[l, shift left=1.8, "i_1^*\circ\Und"] 
		\end{tikzcd}
	\end{center}
	satisfies the conditions of \cite[Proposition 2.2.3]{riehlverity2017inftycosmos}, where $(p_1^*,i_1^*)$ is the adjunction between quasi-categories and complete Segal spaces defined in \cite[Section 4]{joyaltierney2007qcatvssegal}. Indeed, we only need to check that the left adjoint $\Diag\circ p_1^*$ preserves binary products. It suffices to check on generators $\Delta^n$ in $\sSet$. We know that $p_1^*$ preserves binary products \cite[Example 2.2.5]{riehlverity2017inftycosmos}, so we can focus on $\Diag$. However, by definition we have $p_1^*(\Delta^n) = F(n)$ and $\Diag$ preserves binary products of $F(n)$ by assumption. 
	\item \label{item:cosmos} Following the previous item we call the fibrant objects in $\M$ $\infty$-categories and denote them with letters $\C,\D$. 
	Moreover, we will use results about $\infty$-categories as developed in \cite{riehlverity2017inftycosmos,riehlverity2018elements}.
	\item \label{item:map} The fact that $\M$ is an $\infty$-cosmos in particular implies that it is simplicially enriched (but not necessarily a simplicial model category). We denote the enrichment by $\Map$. 
	\item \label{item:und} For a given $\infty$-category $\C$, we call the complete Segal space $\Und(\C)$ the {\it underlying complete Segal space}.
	\item \label{item:discrete} We call an $\infty$-category an {\it $\infty$-groupoid} if it is a discrete $\infty$-category \cite[Definition 1.2.26]{riehlverity2018elements}, which is equivalent to the underlying complete Segal space being a homotopically constant simplicial space.
	\item \label{item:delta} Let $\DDelta$ be the category of simplices, then using the {\it classifying diagram} \cite[Subsection 3.5]{rezk2001css}, we get a complete Segal space $N\DDelta$ and then applying $\Row$ we get an $\infty$-category $\Row N(\DDelta)$. In order to simplify notation, we will simply denote this $\infty$-category by $\DDelta$. Moreover, for a given $\infty$-category $\C$, we denote $\infty$-category of {\it cosimplicial objects} $\kX: \DDelta \to \C$ by $cos\C$. Similarly, we denote the $\infty$-category of {\it simplicial objects} $\kX: \DDelta^{op} \to \C$ by $s\C$. Moreover, we use $\ks,\kt: \kX_1 \to \kX_0$ to denote the {\it source} and {\it target} map. 
	\item \label{item:fn} Let $F(n)$ be the free complete Segal space on $n$ objects \cite[Subsection 2.3]{rezk2001css}. Then we denote the $\infty$-category $\Row(F(n))$ by $\mathbbe{n}$ following notation in \cite{riehlshulman2017rezktypes} \cite[Definition 1.1.4]{riehlverity2018elements}. In particular, $\1$ is the terminal $\infty$-category in $\M$ and we denote the unique map by $!:X \to \1$. Moreover, $\2$ is the ``free arrow". We use $0,1: \1 \to \2$ for the two inclusion maps.
	\item \label{item:fnund} By definition the objects $\mathbbe{n}$ form a cosimplicial object $\DDelta \to \M$ and the resulting functor 
	$\M \to \ss$ that takes $X$ to $\Map_{\ss}(\mathbbe{n},X)$ is (equivalent to) $\Und(X)$. 
	\item For a given $\infty$-category $\C$, an {\it object} is a map $\kx: \1 \to \C$ and a {\it morphism} is a map $\kf: \2 \to \C$. 
	\item \label{item:mapping space} For an $\infty$-category $\C$ and two objects $\kx,\ky$ in $\C$, we define the {\it mapping space} $\map_\C(\kx,\ky)$ as the pullback $\1 \times_{\C \times \C} \C^{\2}$ and observe it is an $\infty$-groupoid \cite[Section 5]{rezk2001css}, \cite[Definition 3.4.9]{riehlverity2018elements}.
	\item \label{item:finitelimits} An $\infty$-category $\C$ has {\it finite limits} if the constant diagram functor has a right adjoint \cite[Definition 2.3.2]{riehlverity2018elements}. This is equivalent to the underlying complete Segal space $\Und(\C)$ having finite limits. 
	If $\C$ has a terminal object then we denote it by $1_\C$.
	\item \label{item:lccc} An $\infty$-category $\C$ with finite limits is {\it locally Cartesian closed} if the pullback functor has a right adjoint \cite[Definition 2.1.1]{riehlverity2018elements}. This is equivalent to $\Und(\C)$  being locally Cartesian closed. If $\C$ is locally Cartesian closed and $\kX$ an object in $\C$, we denote the right adjoint to the product functor $- \times \kX: \C \to \C_{/\kX}$ as $\Gamma_\kX: \C_{/\kX} \to \C$.
	\item \label{item:limits} Let $\C$ be an $(\infty,1)$-category with finite limits. Then {\it $\map_\C(\kD,-)$ reflects limits}, meaning that for every diagram $\kF:I \to \C$ an object $\kL$ in $\C$ is the limit of the diagram if and only if for every object $\kD$ in $\C$, $\map_\C(\kD,\kL)$ is the limit of the diagram of spaces $\map_\C(\kD,\kF)$ \cite[Proposition 4.3.1]{riehlverity2018elements}. 
	\item \label{item:examples} If $\C$ is an established quasi-category (complete Segal space) that appears in the literature (and is obtained by applying the simplicial nerve \cite[Proposition 1.1.5.10]{lurie2009htt} to an important simplicial model structure), then we will denote the corresponding $\infty$-category again by $\C$, instead of $\Row(t^!\C)$ ($\Row(\C)$) to simplify notation. For example, we denote the $\infty$-category of spaces (given as the simplicial nerve of Kan complexes) by $\s$ and the $\infty$-category of small $\infty$-categories (given as the nerve of $\M$ itself) by $\cat_\infty$.
\end{enumerate}
\subsection{Covariant Model Structure} \label{subsec:covar}
Next we review the relevant results regarding the covariant (contravariant) model structure:
\begin{enumerate}
	\item We call the fibrant objects in the covariant (contravariant) model structure on $\mathcal{M}_{/X}$ {\it left fibrations} ({\it right fibrations}).
	\item \label{item:functorial} The covariant model structure is functorial in the base. In particular, for any map of generic objects $X \to Y$ and covariant equivalence $f:A \to B$ over $X$, $f: A \to B$ is also a covariant equivalence over $Y$.
	\item \label{item:left localization} If $\C$ is an $\infty$-category and $\mathscr{F} \to \C$ is a left (right) fibration, then $\mathscr{F}$ is also an $\infty$-category (as the covariant (contravariant) model structure is a Bousfield localization). Notice \cite[Subsection 5.5]{riehlverity2018elements} calls them {\it discrete (co)Cartesian fibrations} of $\infty$-categories.
	\item \label{item:left fib fiber} A map $\mathscr{F} \to \1$ is a left and right fibration if and only if $\mathscr{F}$ is an $\infty$-groupoid (\infref{item:discrete}) \cite[Lemma 5.5.2]{riehlverity2018elements}. Hence, the fiber of each left (right) fibration is an $\infty$-groupoid.
	\item For a generic object $X$ in $\M$ and map $\kX:\1 \to X$ we denote a choice of covariant (contravariant) fibrant replacement by $L_\kX \to X$ ($R_\kX \to X$) and call it the {\it representable left (right) fibration} represented by $\kX$. 
	\item \label{item:undercat} \cite[Corollary 5.5.13]{riehlverity2018elements} If $\C$ is an $\infty$-category and $\kX$ an object, then the under-category $\C_{\kX/}= \C^{\2} \times_\C \1$ is a covariant fibrant replacement via the map $\{\id_\kX\}:\1 \to \C_{\kX/}$ \cite[Theorem 5.7.1]{riehlverity2018elements}. Notice the fiber over an object $\kY$ is the mapping space $\map_\C(\kX,\kY)$. Similarly, we use the over-category $\C_{/\kX}$ as the representable right fibrations and observe it has fiber $\map_\C(\kY,\kX)$ over an object $\kY$. Notice this coincides with the definitions of under-categories for quasi-categories by \cite[Proposition D.6.4]{riehlverity2018elements}.
	\item \label{item:cocone} \cite[Proposition D.6.4]{riehlverity2018elements} For a diagram $\kF: X \to \C$, where $X$ is a generic object and $\C$ is an $\infty$-category, we get a left fibration of cocones $\C_{\kF/} \to \C$ defined as the limit of the diagram 
	\begin{center}
		\begin{tikzcd}[column sep=0.4in]
			\1 \arrow[r, "\kF"] & \C^X & \C^{X \times \2} \arrow[l, "0^*"'] \arrow[r, "1^*"] & \C^X & \C \arrow[l, "!^*"']  
		\end{tikzcd}
	\end{center}
    and the left fibration of cocones $\C_{\kF/}$ is equivalent to a representable left fibration $\C_{\kc/}$ if and only if $\kc$ is the colimit of the diagram $\kF$ \cite[Proposition 4.3.2]{riehlverity2018elements}.
    Similarly, we get a right fibration of cones $\C_{/\kF} \to \C$ defined as the limit of the following diagram 
    \begin{center}
    	\begin{tikzcd}[column sep=0.4in]
    		\C \arrow[r, "!^*"] & \C^X & \C^{X \times \2} \arrow[l, "0^*"'] \arrow[r, "1^*"] & \C^X & \1 \arrow[l, "\kF"']  
    	\end{tikzcd}
    \end{center}
	and the right fibration of cones $\C_{/\kF}$ is equivalent to a representable right fibration $\C_{/\kl}$ if and only if $\kl$ is the limit of the diagram $\kF$.
	\item \label{item:final} \cite[Corollary 3.6.10.]{riehlverity2018elements} Let $\D$ be an $\infty$-category with final object $\kt$. Then for any functor $\kF: \D\to \C$, the inclusion map $\kt: \1 \to \D$ induces an equivalence of left fibrations $\C_{\kF/} \xrightarrow{ \ \simeq \ } \C_{\kF(\kt)/}$. 
\end{enumerate}
	
\subsection{Reedy Covariant Model Structure} \label{subsec:reedy covar}
Finally, we use the covariant (contravariant) structure to define the {\it Reedy covariant (contravariant) model structure} and review some important properties. Let us denote simplicial objects in $\M_{/X}$ and $\ss_{/\Inc(X)}$ by $s\M_{/X}$ and $\sss_{/\Inc(X)}$. Then the adjunction \ref{eq:relative adjunctions} gives us the Quillen equivalences
	\begin{center}
	\adjun{(s\M_{/X})^{cov}}{(\sss_{/\Inc(X)})^{cov}}{s\Inc}{s\Row}, \
	\adjun{(s\M_{/X})^{contra}}{(\sss_{/\Inc(X)})^{contra}}{s\Inc}{s\Row}
\end{center}
where we are giving the categories the Reedy model structure and $s\Row$, $s\Inc$ are define level-wise \cite[Proposition 15.4.1]{hirschhorn2003modelcategories}.

We want to review several important properties of these two model structures, which requires the {\it Segal condition} (\cite[Section 5]{rezk2001css} also reviewed in \cref{def:so}) and the {\it completeness condition} (\cite[Section 6]{rezk2001css}\cite[Section 10]{rezk2010thetanspaces} also reviewed in \cref{def:cso}).
\begin{enumerate}
	\item A fibrant object is called a Reedy left (right) fibration and the weak equivalences Reedy covariant (contravariant) equivalences. 
	\item For a Reedy left fibration $L_\bullet \to X$, we denote by $\LFib_k(L) \to X$  the $k$-th level left fibration.
	Similarly, for a Reedy right fibration $R_\bullet \to X$, we denote by $\RFib_k(R) \to X$ the $k$-th level right fibration.
	\item Generalizing \covref{item:left fib fiber}, for a Reedy left (right) fibration $F_\bullet \to X$, the fiber $\1 \times_X F_\bullet$ is a simplicial $\infty$-groupoid.
	\item \label{item:segal cart} For a given Reedy left (right) fibration $F \to X$ the following are equivalent:
	\begin{enumerate}
		\item The simplicial object $F_\bullet$ satisfies the Segal condition.
		\item The simplicial simplicial space $\Und(F_\bullet)$ satisfies the Segal condition. 
		\item For every object $\kX: \1 \to X$ the simplicial $\infty$-groupoid $\1 \times_X F$ satisfies the Segal condition. 
		\item For every object $\kX: \1 \to X$ the underlying simplicial space $\Und(\1 \times_X F)$ satisfies the Segal condition. 
	\end{enumerate}
	Indeed, the equivalences $(a) \Leftrightarrow (b)$ and $(c) \Leftrightarrow (d)$ follow from the fact that $\Und$ is a Quillen right adjoint of a Quillen equivalence. The equivalence $(b) \Leftrightarrow (d)$ is proven in \cite[Theorem 4.7]{rasekh2021cartfibcss}.
	If these equivalent conditions hold, then we call $F$ a {\it Segal (co)Cartesian fibration}.
	\item \label{item:cart} For a given Segal (co)Cartesian fibration $F_\bullet \to X$ the following are equivalent:
	\begin{enumerate}
		\item The simplicial object $F_\bullet$ satisfies the completeness condition. 
		\item The simplicial simplicial space $\Und(F_\bullet)$ satisfies the completeness condition. 
		\item For every object $\kX: \1 \to X$ the fiber $\1 \times_X F_\bullet$ satisfies the completeness condition. 
		\item For every object $\kX: \1 \to X$ the underlying simplicial space $\Und(\1 \times_X F_\bullet)$ satisfies the completeness condition. 
	\end{enumerate}
	 Again, the equivalences $(a) \Leftrightarrow (b)$ and $(c) \Leftrightarrow (d)$ follow from using $\Und$ and $(b) \Leftrightarrow (d)$ from \cite[Theorem 4.7]{rasekh2021cartfibcss}.
	If these equivalent conditions hold, then we call $F_\bullet$ a {\it (co)Cartesian fibration}.
\end{enumerate}

{\bf Remark regarding Cartesian fibrations:} The definition of coCartesian fibration given here does not directly coincide with the definition in an $\infty$-cosmos given in \cite[Definition 5.2.1]{riehlverity2018elements}, however, they are in fact equivalent.
Indeed, it is proven in \cite[Subsection 3.1]{rasekh2021cartfibmarkedvscso} that the definition in (\reecovref{item:cart}) is equivalent to the definition due to Lurie \cite[Definition 2.4.2.1]{lurie2009htt}, which itself is equivalent to the the definition in an $\infty$-cosmos \cite[Proposition F.4.5]{riehlverity2018elements}. This chain of equivalences allows us to translate any Cartesian fibration and equivalence of Cartesian fibrations as defined here into results about Cartesian fibrations in the $\infty$-cosmos sense.

In particular this means that results, such as the existence of representable Cartesian fibrations (\cref{prop:rep left fib higher cat}/\ref{eq:contra table}) and the {\it Yoneda lemma for representable Cartesian fibrations} (\cref{the:yoneda Lemma for cosimplicial objects}/\ref{eq:contra table}), will also hold for Cartesian fibrations in an $\infty$-cosmos.

\subsection{Examples} \label{subsec:examples}

Let us provide two examples of a theory of $\infty$-categories with fibration structure:
\begin{enumerate}
	\item\label{item:css} We want to prove that complete Segal spaces give us a theory of $\infty$-categories with fibrant structure.
	We will check the three conditions.
	\begin{enumerate}
		\item[(I)] There is a left proper combinatorial model structure compatible with Cartesian closure on the category of simplicial spaces with fibrant objects complete Segal spaces \cite[Theorem 7.2]{rezk2001css}. 
		\item[(II)] We define the functors in \ref{eq:adjunctions} to be the identity functors and all conditions hold evidently.
		\item[(III)] For every simplicial space there exists a simplicial combinatorial left proper model structure on the over-category \cite[Theorem 3.12]{rasekh2017left} that is in fact a left Bousfield localization of the complete Segal space model structure \cite[Theorem 5.11]{rasekh2017left}. Moreover the adjunctions \ref{eq:relative adjunctions} are just identity functors and hence evidently Quillen equivalences.
	\end{enumerate}
    Hence, all conditions hold.
	\item \label{item:qcat} Next, we want to show that quasi-categories give us a theory of $\infty$-categories with fibration structure.
	Again, we check the three conditions.
	\begin{enumerate}
		\item[(I)] There is a combinatorial left proper model structure compatible with Cartesian closure on the category of simplicial sets, the {\it Joyal model structure} with fibrant objects quasi-categories \cite[Theorem 2.2.5.1, Corollary 2.2.5.4]{lurie2009htt}.
		\item[(II)] Next, we define the functors in the diagram \ref{eq:adjunctions} as follows: 
		$(\Diag = t_!, \Und = t^!)$ \cite[Section 2]{joyaltierney2007qcatvssegal}, 
		$(\Inc = p_1^*, \Row = i_1^*)$ \cite[Section 4]{joyaltierney2007qcatvssegal}. Both $(t_!,t^!)$ \cite[Theorem 4.12]{joyaltierney2007qcatvssegal} and $(p_1^*,i_1^*)$ \cite[Theorem 4.11]{joyaltierney2007qcatvssegal} are proven to be Quillen equivalences. Finally $t_!p_1^*$ is the identity functor and so commutes with binary products.
		\item[(III)] Finally, we also have the combinatorial left proper model structures on the category of simplicial sets called the covariant and contravariant model structures \cite[Section 8]{joyal2008theory}\cite[Proposition 2.1.4.7]{lurie2009htt} that are simplicial \cite[Proposition 2.1.4.8]{lurie2009htt} and a Bousfield localization of the Joyal model structure \cite[Theorem 3.1.5.1]{lurie2009htt}. Finally the adjunction $p_1^*,i_1^*$ make both adjunctions in \ref{eq:relative adjunctions} into a Quillen equivalence of covariant (contravariant) model structures \cite[Theorem B.12]{rasekh2017left}.
	\end{enumerate}
	Hence, all conditions hold.
\end{enumerate}

\section{Representable Reedy Left Fibrations} \label{sec:repfib}
In this section we define {\it representable Reedy left fibrations} over generic objects. 
We then give an explicit description of representable Reedy left fibrations over $\infty$-categories, generalizing a similar result for representable left fibrations. Finally, we prove the Yoneda lemma for representable Reedy left fibrations.
We will use the language of an $\infty$-category with fibration structure as introduced in \cref{sec:infty}.

 Let $X$ be a generic object (\infref{item:gen}) and $\kx^\bullet: \DDelta \to X$ be a cosimplicial object (\infref{item:delta}). We want to construct a Reedy left fibration we denote $L_{\kx^\bullet}$ that has the universal property that  $\LFib_k(L_{\kx^\bullet})$ is a covariant fibrant replacement of the map $\{ \kx^k \}: \1 \to X$ (\covref{item:undercat}) and notice this implies this simplicial object is unique up to equivalence. 
 
 Our naive guess might be to simply choose the fibrant replacements level-wise. However, while this certainly gives us a collection of left fibrations it would not give us a Reedy left fibration as the various left fibrations do not interact in any way. We need a more global approach, meaning we need an appropriate analogue of the point $\1$ that we can use to construct fibrant replacements.
 
 Let $\pi^\bullet_i: \DDelta_{\bullet /} \to \DDelta$ be the Reedy left fibration that is defined level-wise as $\LFib_k(\pi^\bullet_i) = \DDelta_{[k]/}$.
 and for a map $\delta: [m] \to [n]$, define $\delta^*: \DDelta_{[n]/} \to \DDelta_{[m]/}$ via precomposition. This Reedy left fibration has the interesting property that it is level-wise representable, which justifies following definition and it subsequent lemma.
 
 \begin{defone} \label{def:rep reedy left fib}
 	A Reedy left fibration $L \to X$ is called {\it representable} if there exists a cosimplicial object $\kx^\bullet: \DDelta \to X$ and a Reedy covariant equivalence $f: \DDelta_{\bullet/} \to L$ over $X$.
 \end{defone}

 \begin{lemone}
 	Let $L \to X$ be a representable Reedy left fibration represented by $\kx^\bullet: \DDelta \to X$. Then for each $k \geq 0$ we have an equivalence of left fibrations $\LFib_k(L) \simeq L_{\kx^k}$.
 \end{lemone}

 \begin{proof}
 	By \cref{def:rep reedy left fib}, we have a covariant equivalence $\DDelta_{[k]/} \to \LFib_k(L)$ over $X$. Moreover, by definition $\1 \to \DDelta_{[k]/}$ is a covariant equivalence over $X$ (\covref{item:undercat}). Hence, $\1 \to \LFib_k(L)$ is a covariant equivalence over $X$. 
 	The result now follows from the definition of fibrant replacements.
 \end{proof}

 We now want to focus on the case of $\infty$-categories (\infref{item:cosmos}) and give an explicit description. As we stated before for an $\infty$-category $\C$ and object $c$, the under-category $\C_{c/}$ is the representable left fibration (\covref{item:undercat}). We want to generalize this observation to representable Reedy left fibrations. However, we again face the problem that for a given cosimplicial object $\kx: \DDelta \to \C$, the various over-categories $\C_{\kx^k/}$ do not form a simplicial object. We need to expand our diagram again. 
 
 \begin{remone} \label{rem:segal also works}
  If we choose our model of $\infty$-categories with fibration structure to be complete Segal spaces (\cref{subsec:examples}), then the results \cref{prop:rep left fib higher cat}, \cref{the:yoneda Lemma for cosimplicial objects} and \cref{cor:equiv of simp obj} would also hold for Segal spaces, which is a weaker requirement than being an $\infty$-category (in this context complete Segal spaces). 
  Indeed, this follows from the fact that under-categories are already the representable left fibrations for Segal spaces as proven in \cite[Theorem 3.49]{rasekh2017left}. 
 \end{remone}

 Let $\pi^\bullet_f: \DDelta_{/\bullet} \to \DDelta$ be the cosimplicial $\infty$-category that is defined level-wise as $\DDelta_{/[k]} \to \DDelta$ and for a map $\delta: [m] \to [n]$, define $\delta^*: \DDelta_{/[n]} \to \DDelta_{/[m]}$ via post-composition. For an $\infty$-category $\C$ and a cosimplicial object $\kx: \DDelta \to \C$, we can use the cosimplicial object $\pi^\bullet_f$ to define the Reedy left fibration 
 $\C_{\kx^\bullet/} = \C_{\pi^\bullet_f \kx^\bullet/}$. We now have following result.
 
 \begin{propone} \label{prop:rep left fib higher cat}
 	Let $\C$ be an $\infty$-category and $\kx^\bullet: \DDelta \to \C$ a cosimplicial object. 
 	Then the map $\C_{\pi_f\kx^\bullet/} \to \C$ is the representable Reedy left fibration represented by $\kx^\bullet$.
 \end{propone}

\begin{proof}
	Based on \covref{item:undercat} the map $\C_{\pi_f^k \kx^k/} \to \C$ is a left fibration for each $k$ and so the map $\pi_f^{\bullet}$ is a Reedy left fibration.
	We will now construct a map of simplicial objects $\DDelta_{\bullet/} \to \C_{\pi_f\kx^\bullet/}$ over $\C$ and prove it is a level-wise covariant equivalence and this this prove that $\C_{\pi_f^k \kx^k/} \to \C$ is in fact the representable Reedy left fibration represented by $\kx^\bullet$.
	
	Define the functor of simplicial categories 
	$$\comp: \DDelta_{\bullet /} \to {(\DDelta^\2)}^{\DDelta_{/\bullet}}$$
	given by $\comp(\delta: [k] \to [m])(\gamma: [n] \to [k]) = \delta \circ \gamma: [n] \to [m]$.
	We can restrict $\comp$ to a functor of simplicial categories
	$$\comp: \DDelta_{\bullet/} \to \DDelta_{\pi^\bullet_f/}. $$
	Indeed, if we restrict the image of $\comp([k] \to [m])$ to the domain in $\DDelta^\2$ we get the projection map $\pi^\bullet_f:\DDelta_{/\bullet} \to \DDelta$. 
	Similarly, if we restrict the image of $\comp$ to the codomain we get a constant map $\DDelta_{/\bullet} \to \DDelta$ with value $[m]$.
	Finally, we get the desired map $\DDelta_{\bullet /}\to\C_{\pi^k_f\kx^k/}$ by post-composing the map $\comp$ with the restricted codomain with the cosimplicial object $\kx^\bullet: \DDelta \to \C$.
	
	We will now proceed to prove this map is a Reedy covariant equivalence. It suffices to prove we have a level-wise equivalence. 
	For each $k$, we get following diagram over $\C$.
	\begin{center}
		\begin{tikzcd}[row sep=0.4in, column sep=0.6in]
			\1 \arrow[drr, "x_k"'] \arrow[r, "\id_k" near end, "\simeq"'] \arrow[rrrr, "\id_{\kx_k}" near start, "\simeq" near end, bend left=20] &  
			\DDelta_{k/}  \arrow[rr, "\kx^k \circ \comp^k"] \arrow[dr] & &
			\C_{\pi^k_f\kx^k/} \arrow[r, "\simeq"] \arrow[dl, twoheadrightarrow] & 
			\C_{\kx^k/} \arrow[dll, twoheadrightarrow] \\
			& & \C &
		\end{tikzcd}
	\end{center}
	The map $\1 \to \DDelta_{k/}$ is a covariant equivalence over $\DDelta$ (by \covref{item:undercat}) and so also a covariant equivalence over $\C$ (\covref{item:functorial}). The map $\1 \to \C_{\kx^k/}$ is directly a covariant equivalence over $\C$ (by \covref{item:undercat}).
	The map of left fibrations $\C_{\pi^k_f\kx^k/}\to \C_{\kx^k/}$ is an equivalence because $\id_{\kx}$ is the terminal object in the diagram (\covref{item:final}). Thus, by $2$-out-of-$3$, the middle map is also a covariant equivalence over $\C$.
\end{proof}

\begin{notone} \label{not:fiber mapping}
 Given that the fiber of $\C_{\pi_f\kx^\bullet/} \to \C$ over an object $\ky$ is level-wise equivalent to $\map_\C(\kx^k,\ky)$ (\covref{item:undercat}), we will henceforth denote this fiber by $\map_\C(\kx^\bullet,\ky)$ or $\map_\C(\kx,\ky)$. 
\end{notone}

Finally, we want to prove the Yoneda lemma for representable Reedy left fibrations. We need following lemma.

\begin{lemone} \label{lemma:pushout of cats}
	The following square is a pushout square of categories
	\begin{center}
		\begin{tikzcd}
			\partial \2 \times \DDelta \arrow[r] \arrow[d] & \partial \2 \times \DDelta \arrow[d] \\
			\2 \times \DDelta_{[k]/} \arrow[r] & \2 \times \DDelta  
		\end{tikzcd}
	\end{center}
\end{lemone}

\begin{proof}
	For a category $\D$ a diagram of the form 
		\begin{center}
		\begin{tikzcd}
			\partial \2 \times \DDelta \arrow[r] \arrow[d] & \partial \2 \times \DDelta \arrow[d] \arrow[ddr, "F", bend left=20] & \\
			\2 \times \DDelta_{[k]/} \arrow[r] \arrow[drr, "G", bend right=20]& \2 \times \DDelta  \arrow[dr, dashed, "\exists ! H" description] & \\
			& & \D 
		\end{tikzcd}
	\end{center}
	has a unique lift. Indeed the fact that the diagram commutes implies that for two objects $[k] \to [m]$, $[k] \to [n]$ and $i = 0,1$ we have 
	$G(i,[k] \to [m]) = G(i,[k] \to [n])$ and so $G$ factors through $\2 \times \DDelta$.
\end{proof}

\begin{theone} \label{the:yoneda Lemma for cosimplicial objects}
	Let $\C$ be an $\infty$-category and $\kx^{\bullet}$ and $\ky^{\bullet}$ be two cosimplicial objects, then we have an equivalence
	$$ \Map_{/\C}(\C_{\pi^{\bullet}_f\kx^{\bullet}/},\C_{\pi^{\bullet}_f\ky^{\bullet}/}) \overset{\simeq}{\twoheadrightarrow} \map_{cos\C}(\ky^{\bullet},\kx^{\bullet})$$
\end{theone}

\begin{proof}
   First, using \cref{prop:rep left fib higher cat} and the fact that $\C_{\pi^{\bullet}_f\ky^{\bullet}/} \to \C$ is a Reedy left fibration we have a trivial Kan fibration 
   \begin{equation} \label{eq:triv kan fib}
   	\Map_{/\C}(\C_{\pi^{\bullet}_f\kx^{\bullet}/},\C_{\pi^{\bullet}_f\ky^{\bullet}/}) \overset{\simeq}{\twoheadrightarrow} \Map_{/\C}(\DDelta_{\bullet/},\C_{\pi^{\bullet}_f\ky^{\bullet}/})
   \end{equation}
   which means we need to understand the space on the right hand side. Before we can do so we need to make a technical construction.
   
   Define $\Cone(\C)$ to be the subspace of functors $ F:\DDelta_{\bullet /} \times \DDelta_{\bullet /} \times \2 \to \C$ such that 
   $H$ fits into a diagram of the following shape:
   \begin{equation} \label{eq:diagram}
   	\begin{tikzcd}[row sep=0.25in, column sep=0.3in]
   		\DDelta_{\bullet /} \times \DDelta_{\bullet /} \arrow[d, hookrightarrow] \arrow[r, "\pi_1"] & \DDelta \arrow[dr, "F"]  & \\
   		\2 \times \DDelta_{\bullet /} \times \DDelta_{\bullet /} \arrow[rr, "H" description] & & \C \\
   		\DDelta_{\bullet /} \times \DDelta_{\bullet /} \arrow[u, hookrightarrow] \arrow[r, "\pi_1"] & \DDelta \arrow[ur, "G"']  &
   	\end{tikzcd}
   \end{equation}
   Notice a point in $\Cone(\C)$ can be characterized as a map out of the pushout of the diagram 
   \begin{center}
   	\begin{tikzcd}[column sep=0.3in]
   		\DDelta &
   		\DDelta_{\bullet/} \times \DDelta_{\bullet /} \arrow[l, "\pi_1"'] \arrow[r, hookrightarrow, "0 \times \id"] &
   		\2 \times \DDelta_{\bullet /} \times \DDelta_{\bullet /} \arrow[r, hookleftarrow, "1 \times \id"] & 
   		\DDelta_{\bullet/} \times \DDelta_{\bullet /} \arrow[r, "\pi_1"] &
   		\DDelta	
   	\end{tikzcd}
   \end{center}
   This means that in order to better understand $\Cone(\C)$ we need to compute this colimit. 
   We can reformulate it as the following two pushout squares and compute it directly
   \begin{center}
   	\begin{tikzcd} 
   		\partial \2 \times \DDelta_{\bullet/} \times \DDelta_{\bullet /} \arrow[d] \arrow[r]  & 
   		\partial \2 \times \DDelta_{\bullet /} \arrow[r] \arrow[d] & 
   		\partial \2 \times \DDelta \arrow[d] \\  
   		 \2 \times \DDelta_{\bullet/} \times \DDelta_{\bullet /} \arrow[r] & 
   		 \2 \times \DDelta_{\bullet/} \arrow[r]  \arrow[ul, phantom, "\ulcorner", very near start] & 
   		 \2 \times \DDelta \arrow[ul, phantom, "\ulcorner", very near start]
   	\end{tikzcd}
   \end{center}
	where the second pushout square follows from applying \cref{lemma:pushout of cats} level-wise.
	Hence, we have a bijection of spaces $\Cone(\C) \cong \Map(\2 \times \DDelta, \C)$.
	We can extend this bijection to a commutative triangle
	\begin{center}
		\begin{tikzcd}
			\Cone(\C) \arrow[rr, "\cong"] \arrow[dr] & & \Map(\2 \times \DDelta, \C) \arrow[dl] \\
			 & \Map(\DDelta,\C) \times \Map(\DDelta,\C) & 
		\end{tikzcd}
	\end{center}
	where both projections are induced by restricting along the endpoints of $\2$. Taking the fiber of the bijection over the points $(\ky,\kx)$ in $ \Map(\DDelta,\C) \times \Map(\DDelta,\C)$, we get a bijection 
	$$\Cone(\C)_{\ky,\kx} \cong \map_{cos\C}(\ky,\kx)$$
	where the left hand side by definition denotes the fiber and the right hand side follows from the definition of the mapping space \infref{item:mapping space}.
	
	Hence, in order to finish the proof it suffices to prove that the fiber $\Cone(\C)_{\ky,\kx}$ is bijective to the space   $\Map_{/\C}(\DDelta_{\bullet/},\C_{\pi^{\bullet}_f\ky^{\bullet}/})$ given in Equation \ref{eq:triv kan fib}.
	Using the fact that $\Map(\DDelta_{\bullet/},-)$ commutes with limits and and the limit description of $\C_{\pi^{\bullet}_f\ky^{\bullet}/}$ in \covref{item:cocone}, we can describe $\Map_{/\C}(\DDelta_{\bullet/},\C_{\pi^{\bullet}_f\ky^{\bullet}/})$ as the following limit diagram 
	\begin{center}
		\begin{tikzcd}
			\Map(\DDelta_{\bullet /},\1) \arrow[dr, "\ky"]  & & \Map(\2 \times \DDelta_{\bullet /},  \C^{\DDelta_{\bullet /}}) \arrow[dl] \arrow[dr] & & \Map(\DDelta_{\bullet /},\1) \arrow[dl, "\kx"] \\
			& \Map(\DDelta_{\bullet /} ,\C^{\DDelta_{\bullet /}}) & & \Map(\DDelta_{\bullet /},\C^{\DDelta_{\bullet /}})
		\end{tikzcd}.
	\end{center}
	Now, using the fact that $\1$ is the final object and the Cartesian closure, we get following limit diagram of spaces:
	\begin{center}
		\begin{tikzcd}
			\1 \arrow[dr, "\ky"]  & & \Map(\2 \times \DDelta_{\bullet /} \times \DDelta_{\bullet /},  \C) \arrow[dl] \arrow[dr] & & \1 \arrow[dl, "\kx"] \\
			& \Map(\DDelta_{\bullet /}\times \DDelta_{\bullet /},\C) & & \Map(\DDelta_{\bullet /}\times \DDelta_{\bullet /},\C)
		\end{tikzcd}.
	\end{center}
	Hence, $\Map_{/\C}(\DDelta_{\bullet/},\C_{\pi^{\bullet}_f\ky^{\bullet}/})$ is precisely the sub-space of $\Map(\2 \times \DDelta_{\bullet /} \times \DDelta_{\bullet /}, \C)$ consisting of diagrams of the form given in \ref{eq:diagram}, such that $F= \ky$ and $G = \kx$, which is precisely $\Cone(\C)_{\ky,\kx}$. This finishes the proof.
\end{proof}

This theorem also has a valuable corollary.

\begin{corone} \label{cor:equiv of simp obj}
	Let $\C$ be an $\infty$-category and $\kx,\ky: \DDelta \to \C$ be two cosimplicial objects. Then $\kx,\ky$ are equivalent in $cos\C$ if and only if the corresponding representable Reedy left fibrations $\C_{\kx/}$ and $\C_{\ky/}$ are equivalent.
\end{corone}

Notice, if $\kx^\bullet, \ky^\bullet$ are discrete simplicial objects in $\C$, then \cref{the:yoneda Lemma for cosimplicial objects} recovers the classical Yoneda lemma of $\infty$-categories \cite[Corollary 5.7.16]{riehlverity2018elements}.
We end this section by observing that all these results also hold in a contravariant setting. We simply use following table of conversions:

\begin{equation} \label{eq:contra table}
	\begin{tabular}{|c|c|}
		\hline 
		Reedy Left Fibration & Reedy Right Fibration \\ \hline 
		cosimplicial object cosX & simplicial object sX \\ 
		$\kx^\bullet: \DDelta \to X$ & $\kx_\bullet: \DDelta^{op} \to X$ \\ \hline
		$L_{\kx^\bullet}$ Representable & $R_{\kx_\bullet}$ Representable\\ 
		$\LFib_k(L_{\kx^\bullet}) \simeq L_{\kx^k}$ & $\RFib_k(R_{\kx_\bullet}) \simeq R_{\kx_k}$ \\ \hline 
		$\pi_i^\bullet:  \DDelta_{\bullet /} \to  \DDelta$ & $\pi^i_\bullet: (\DDelta^{op})_{/\bullet } \to  \DDelta^{op}$ \\ \hline
		$\pi_f^\bullet: \DDelta_{/\bullet} \to \DDelta$ & $\pi^f_\bullet: (\DDelta^{op})_{\bullet /} \to \DDelta^{op}$ \\ \hline 
		Over $\infty$-Categories: & Over $\infty$-Categories: \\ 
		 $\C_{\pi^{\bullet}_f\kx^{\bullet}/} \to \C$ & $\C_{/\pi_{\bullet}^f\kx_{\bullet}} \to \C$ \\ \hline 
		$ \Map_{/\C}(\C_{\pi^{\bullet}_f\kx^{\bullet}/},\C_{\pi^{\bullet}_f\ky^{\bullet}/}) \overset{\simeq}{\twoheadrightarrow} \map_{cos\C}(\ky^{\bullet},\kx^{\bullet})$ & 
		$ \Map_{/\C}(\C_{/\pi_{\bullet}^f\kx_{\bullet}},\C_{/\pi_{\bullet}^f\ky_{\bullet}}) \overset{\simeq}{\twoheadrightarrow} \map_{s\C}(\kx_{\bullet},\ky_{\bullet})$ \\ 
		$\C_{\pi^{\bullet}_f\kx^{\bullet}/} \simeq \C_{\pi^{\bullet}_f\ky^{\bullet}/} \Leftrightarrow \kx^\bullet \simeq \ky^\bullet$ & 
		$\C_{/\pi_{\bullet}^f\kx_{\bullet}} \simeq \C_{/\pi_{\bullet}^f\ky_{\bullet}} \Leftrightarrow \kx_\bullet \simeq \ky_\bullet$ \\ \hline 
	\end{tabular}    
\end{equation}

\section{Complete Segal Objects and Representable Cartesian Fibrations}
In this section we study a model of internal $\infty$-categories, {\it complete Segal objects}, using representable Cartesian fibrations. 
In particular, we develop the category theory of {\it complete Segal objects} generalizing the category theory of complete Segal spaces as developed in \cite{rezk2001css}. Throughout this section $\C$ is an $\infty$-category (\infref{item:cosmos}) with finite limits (\infref{item:finitelimits}).

\begin{defone} \label{def:so}
	Let $\kW: \DDelta^{op} \to \C$ be a simplicial object in $\C$. Then $\kW$ is a {\it Segal object} if the map
	$$\kW_n \to \kW_1 \times_{\kW_0} ... \times_{\kW_0} \kW_1$$
	is an equivalence in $\C$ \cite[Section 5]{rezk2001css}, meaning $\kW_n$ is the limit of the diagram 
	\begin{center}
		\begin{tikzcd}
			\kW_1 \arrow[r, "\kt"] & \kW_0 & \kW_1 \arrow[l, "\ks"'] \arrow[r] & ... & \kW_1 \arrow[l] \arrow[r, "\kt"]& \kW_0 & \kW_1 \arrow[l, "\ks"']
		\end{tikzcd}
	\end{center}
	where $\ks,\kt:\kW_1 \to \kW_0$ are the source and target map \infref{item:delta}.
\end{defone}

We now want to study the category theory of Segal objects following the definitions in \cite[Section 5]{rezk2001css}
An {\it object} in $\kW$ is a morphism $\kx: \kD \to \kW_0$. We say $\kD$ is the {\it context} of the object $\kx$.
A {\it morphism} in $\kW$ is a morphism $\kf: \kD \to \kW_1$ in $\C$.

Next we want to define mapping objects, however, for that we need to assume that $\C$ is locally Cartesian closed. With this assumption for two objects $\kx,\ky:\kD \to \kW$ with context $\kD$, define the {\it mapping object} $\map_\kW(\kx,\ky)$ in $\C$ as follows
$$\map_\kW(\kx,\ky) = \Gamma_\kD (\kx,\ky)^* \kW_1 $$
where $\Gamma_\kD$ is the global section functor \infref{item:lccc}.
Notice, by definition we have an equivalence 
$$
\map_\C(1_\C,\map_\kW(\kx,\ky)) \simeq \map_{/\kD}(\id_\kD, (\kx,\ky)^*(\ks,\kt)) \simeq \map_{/\kW_0 \times \kW_0}(\kD,\kW_1)
$$
where $1_\C$ is the terminal object (\infref{item:finitelimits}) and $\ks,\kt:\kW_1 \to \kW_0$ are the source target maps (\infref{item:delta}).
Hence a point in $\map_\kW(\kx,\ky)$ is precisely a map $\kD \to \kW_1$ over $(\kx,\ky): \kD \to \kW_0 \times \kW_0$, as we wanted. 

We want to proceed to study homotopy equivalences in a Segal object. Let $\kZ_\kW(3)$ be the limit of the diagram $Z(3) \to \C$ of the following shape
\begin{center}
	\begin{tikzcd}
		\kW_1 \arrow[r, "\ks"] & \kW_0 & \kW_1 \arrow[l, "\ks"'] \arrow[r, "\kt"] & \kW_0 & \kW_1 \arrow[l,"\kt"'] 
	\end{tikzcd}
\end{center}
and notice it comes with maps $\kW_1 \to \kZ_\kW(3)$ and $\kW_3 \to \kZ_W(3)$. 
Define  $\kW_{hoequiv}$ as the pullback 
\begin{equation}\label{eq:hoequiv}
	\begin{tikzcd}[row sep=0.4in, column sep=0.5in]
		\kW_{hoequiv} \arrow[r] \arrow[d, "\ki"] \arrow[dr, phantom, "\ulcorner", very near start]& \kW_3 \arrow[d, "(d_1d_3 \comma d_0d_3 \comma d_1d_0)"] \\
		\kW_1 \arrow[r, "(s_0d_1\comma \id \comma s_0d_0)"] & \kZ_\kW(3)
	\end{tikzcd}
\end{equation}
We say a morphism $\kf: \kD \to \kW_1$ is a {\it homotopy equivalence} in $\kW$ if the map factors through the map $\kW_{hoequiv}$, which is equivalent to saying that $\map_{/\kW_1}(\kf: \kD \to \kW_1, i: \kW_{hoequiv} \to \kW_1)$ is non-empty.
Finally, notice the unique maps $\kW_0 \to \kW_1, \kW_0 \to \kW_3$ make Diagram \ref{eq:hoequiv} into a commutative square and so we have a map $\kj: \kW_0 \to \kW_{hoequiv}$.

\begin{defone} \label{def:cso}
	A Segal object $\kW$ is a {\it complete Segal object} if the map $\kj:\kW_0 \to \kW_{hoequiv}$ is an equivalence in $\C$ \cite[Section 6]{rezk2001css}.
\end{defone}

We now want to proceed and prove results about complete Segal objects that we know hold for complete Segal spaces. 
For example, if $W$ is a Segal space, then the map of spaces  $i: W_{hoequiv} \to W_1$ is in fact an inclusion of path-components \cite[Lemma 5.8]{rezk2001css}. More importantly, we can effectively characterize equivalences of complete Segal spaces as the fully faithful essentially surjective functors \cite[Proposition 7.6]{rezk2001css}. In order to generalize such results we want to study the relation between (complete) Segal object to Reedy right fibration it represents, using the material from \cref{sec:repfib}.

\begin{propone} \label{prop:segal obj vs Segal Cart}
	Let $\kW: \DDelta^{op} \to \C$ be a simplicial object. 
	Then $\kW$ is a (complete) Segal object if and only if $\C_{/\kW} \to \C$ is a (Segal) Cartesian fibration.
\end{propone}

\begin{proof}
	We already know that $\C_{/\kW}$ is a Reedy right fibration (\cref{prop:rep left fib higher cat}/\ref{eq:contra table}).
	By \reecovref{item:segal cart}/\reecovref{item:cart}, $\C_{/\kW}$ is a (Segal) Cartesian fibration if and only if $\map_\C(\kD,\kW)$ is a (complete) Segal space for all $\kD$. 
	
	The result now follows from the fact that $\map_\C(\kD,-)$ reflects limits \infref{item:limits} and the fact that both conditions of a (complete) Segal object are limits conditions (\cref{def:so},\cref{def:cso}). 
\end{proof}
	Combining this result with \cref{cor:equiv of simp obj}/\ref{eq:contra table}, which implies that two complete Segal objects are equivalent if and only if their corresponding representable Cartesian fibrations are equivalent, means that we can use Cartesian fibrations to study complete Segal objects. 
	In particular, we can translate the category theory of complete Segal objects we developed until now into the language of Cartesian fibrations.
	
	Let $\C$ be an $\infty$-category and let $\kW$ be a Segal object in $\C$, with representable Segal Cartesian fibration $\C_{/\kW} \to \C$. 
	An object $\kx$ in $\kW$ with context $\kD$, $\kx: \kD \to \kW_0$, corresponds to an object in the complete Segal space $\map_\C(\kD,\kW)$. 
	Similarly, a morphism in $\kW$ with context $\kD$, $\kf: \kD \to \kW_1$ corresponds to a morphism in the complete Segal space $\map_\C(\kD,\kW)$.
	Finally, for two objects $\kx,\ky: \kD \to \kW_0$, the mapping object $\map_\kW(\kx,\ky)$ satisfies the universal property 
  	 \begin{equation} \label{eq:mapping object to spaces}
  	 	\map_{\C}(1_\C,map_{\kW}(\kx,\ky)) \simeq  \map_{\map_\C(\kD,\kW)}(\kx,\ky)
  	 \end{equation}
  	meaning that a point in $\map_\kW(\kx,\ky)$ corresponds precisely to a morphism in the Segal space $\map_\C(\kD,\kW)$ with domain $\kx$ and codomain $\ky$. Finally, we also have following relation of homotopy equivalences 

\begin{lemone}\label{lemma:equiv}
	Let $\kW$ be a Segal object. Then for every object $\kD$, there is an equivalence of spaces $\map_\C(\kD,\kW_{hoequiv}) \simeq \map_\C(\kD,\kW)_{hoequiv}$.
	Thus a morphism $\kf:\kD \to \kW_1$ is a homotopy equivalence in $\kW$ if and only if $\kf$ is a homotopy equivalence in the complete Segal space $\map_\C(\kD,\kW)$.
\end{lemone}

\begin{proof}
	This follows directly from the fact that $\map_\C(\kD,-)$ reflects limits (\infref{item:limits}) combined with the fact that $\kW_{hoequiv}$ is characterized via a limit diagram (\ref{eq:hoequiv}).
\end{proof}

     We can now also prove the desired result about homotopy equivalences.
     
\begin{propone} \label{prop:hoequivintomap}
	Let $\kW$ be a Segal object in $\C$. Then the map $\ki:\kW_{hoequiv} \to \kW_1$ is an inclusion.
\end{propone}

\begin{proof}
	The map is an inclusion if and only if the following is a pullback square in $\C$
	\begin{center}
		\begin{tikzcd}
			\kW_{hoequiv} \arrow[r, "="] \arrow[d, "="] \arrow[dr, phantom, "\ulcorner", very near start] & \kW_{hoequiv} \arrow[d, "\ki"] \\
			\kW_{hoequiv} \arrow[r, "\ki"] & \kW_1
		\end{tikzcd}.
	\end{center}
	The result now follows from the fact that $\map_\C(\kD,-)$ reflects limits (\infref{item:limits}) and the the fact that $\map_\C(\kD,\kW)$ is a Segal space (\cref{prop:segal obj vs Segal Cart}). 
\end{proof}

We move on to the final goal of this section: Characterizing the equivalences of complete Segal objects and prove they are analogous to the Dwyer-Kan equivalences of complete Segal spaces \cite[Proposition 7.6]{rezk2001css}.
This result has been called {\it ``fundamental theorem of quasi-categories" } \cite{rezk2017qcats}, so accordingly we call it the 
{\it fundamental theorem of complete Segal objects}.

Let $\C$ be an $\infty$-category with finite limits. A functor of complete Segal objects $\kF: \kW \to \kV$ is {\it fully faithful} if the square
\begin{equation}\label{eq:ff}
	\begin{tikzcd}
		\kW_1 \arrow[r] \arrow[d] \arrow[dr, phantom, "\ulcorner", very near start] & \kV_1 \arrow[d] \\
		\kW_0 \times \kW_0 \arrow[r] & \kV_0 \times \kV_0 
	\end{tikzcd}
\end{equation} 
is a pullback square. Moreover, it is {\it essentially surjective} if for any object $\ky: \kD \to \kV_0$ with context $\kD$, there exists an object $\kx: \kD \to \kW_0$ with context $\kD$ such that $\kF\kx$ is equivalent to $\ky$ in $\kV$.

These definitions translate appropriately to fibrations.

\begin{lemone} \label{lemma:ff}
	Let $\kF:\kW \to \kV$ be a functor of complete Segal objects. Then $\kF$ is fully faithful if and only if for every objects $\kD$, $\map_\C(\kD,\kF): \map_\C(\kD,\kW) \to \map_\C(\kD,\kV)$ is a fully faithful map of complete Segal spaces.
\end{lemone}

\begin{proof}
 Combining the fact that $\map(\kD,-)$ reflects limits (\infref{item:limits}) and the  pullback square \ref{eq:ff} implies that we have a homotopy pullback square of spaces 
 \begin{center}
 	\begin{tikzcd}
 		\map_\C(\kD,\kW)_1 \arrow[r] \arrow[d] \arrow[dr, phantom, "\ulcorner", very near start] & \map_\C(\kD,\kV)_1 \arrow[d] \\
 		\map_\C(\kD,\kW)_0 \times \map_\C(\kD,\kW)_0 \arrow[r] & \map_\C(\kD,\kV)_0 \times \map_\C(\kD,\kV)_0 
 	\end{tikzcd}.
 \end{center}
 Hence, it suffices to prove that this pullback square is equivalent to $\map(\kD,\kF): \map(\kD,\kW) \to \map(\kD,\kV)$ being a fully faithful. 
 We can rephrase this homotopy pullback square as stating that the horizontal map 
 \begin{center}
  \begin{tikzcd}[column sep=-0.5in]
  	\map_\C(\kD,\kW)_1 \arrow[rr, "\simeq"] \arrow[dr, twoheadrightarrow] & & (\map_\C(\kD,\kW)_0 \times \map_\C(\kD,\kW)_0) \times_{\map_\C(\kD,\kV)_0 \times \map_\C(\kD,\kV)_0 }\map_\C(\kD,\kV)_1 \arrow[dl, twoheadrightarrow]  \\
  	& \map_\C(\kD,\kW)_0 \times \map_\C(\kD,\kW)_0
  \end{tikzcd}
 \end{center}
 is an equivalence. This is equivalent to being a fiber-wise equivalence over $\map_\C(\kD,\kW)_0 \times \map_\C(\kD,\kW)_0$. However, for a given point $(\kx,\ky)$ in $\map_\C(\kD,\kW)_0 \times \map_\C(\kD,\kW)_0$, the fiber is precisely $\map_{\map_\C(\kD,\kW)}(\kx,\ky) \to \map_{\map_\C(\kD,\kV)}(\kF\kx,\kF\ky)$ giving us the more common definition of being a fully faithful functor of complete Segal spaces.
\end{proof}

\begin{lemone}\label{lemma:es}
	Let $\kF:\kW \to \kV$ be a functor of complete Segal objects. Then $\kF$ is essentially surjective if and only if for every objects $\kD$, $\map_\C(\kD,\kF): \map_\C(\kD,\kW) \to \map_\C(\kD,\kV)$ is an essentially surjective functor of complete Segal spaces.
\end{lemone}

\begin{proof}
 By \cref{lemma:equiv} weak equivalences in the complete Segal object $\kW$ and the complete Segal space $\map_\C(\kD,\kW)$ agree with each other. 
 Hence, the results follows. 
\end{proof}

We can use our characterizations via complete Segal spaces to give an alternative characterization of fully faithful functors in the locally Cartesian closed setting.

\begin{lemone}
	Let $\C$ be locally Cartesian closed and $\kF: \kW \to \kV$ be a functor of complete Segal objects. 
	Then $\kF$ is fully faithful if and only if for every two objects $\kx,\ky: \kD \to \kW$ in $\kW$, the induced map
	\begin{equation}\label{eq:some map}
		\map_\kW(\kx,\ky) \to \map_\kV(\kF\kx,\kF\ky)
	\end{equation}
	is an equivalence. 
\end{lemone}

\begin{proof}
	First assume  the maps in \ref{eq:some map} are equivalences in $\C$, then by the equivalence \ref{eq:mapping object to spaces} we get an equivalence of mapping spaces of complete Segal spaces for every object $\kD$. Hence, by \cref{lemma:ff}, $\kF$ is fully faithful.
	
	On the other side, assume that $\kF$ is fully faithful. Then for any object two objects $\kx,\ky: \kD \to \kW_0$ we have following pullback squares in $\C$
	\begin{center}
		\begin{tikzcd}
			(\kx \comma \ky)^*\kW_1 \arrow[r] \arrow[d] \arrow[dr, phantom, "\ulcorner", very near start]  & \kW_1 \arrow[r, "\kF_1"] \arrow[d] \arrow[dr, phantom, "\ulcorner", very near start] & \kV_1 \arrow[d] \\
			\kD \arrow[r, "(\kx \comma \ky)"] & \kW_0 \times \kW_0 \arrow[r, "\kF_0 \times \kF_0"] & \kV_0 \times \kV_0
		\end{tikzcd}.
	\end{center}
	By the pasting lemma the rectangle is a pullback as well. Hence we get an equivalence 
	$$ (\kx,\ky)^*\kW_1 \to (\kF\kx,\kF\ky)^*\kV_1$$
	which implies that \ref{eq:some map} is an equivalence as well.
\end{proof}

\begin{theone} \label{the:fundamental theorem cso}
	Let $\kF: \kW \to \kV$ be a functor of complete Segal objects. Then the following are equivalent.
	\begin{enumerate}
		\item $\kF$ is an equivalence of complete Segal objects.
		\item $\kF$ is fully faithful and essentially surjective.
		\item For every object $\kD$, the induced map $\map_\C(\kD,\kF)$ is an equivalence of complete Segal spaces.
		\item For every object $\kD$, the induced map $\map_\C(\kD,\kF)$ is fully faithful and essentially surjective.
	\end{enumerate}
\end{theone}

\begin{proof}
	(2) $\Leftrightarrow$ (4) 
	Follows from combining \cref{lemma:ff} and \cref{lemma:es}.
	
	(1) $\Leftrightarrow$ (3) Is the statement of \cref{cor:equiv of simp obj}
	
	(3) $\Leftrightarrow$ (4) Due to Rezk \cite[Proposition 7.6]{rezk2001css}.
\end{proof}

	From the perspective of this proof we can see why it would not have sufficed to consider an object in a complete Segal object 
	$\kW$ to be a map out of the final object $1_\C \to \kW_0$. We do need to consider objects with different contexts to be able to understand equivalences of complete Segal object in terms of the corresponding complete Segal spaces.

 We end this section by looking at some examples of complete Segal objects in various $\infty$-categories. 
 \begin{enumerate}
 	\item Let $\s$ be the $\infty$-category of spaces (\infref{item:examples}). Then a (complete) Segal object in $\s$ corresponds to a (complete) Segal space (up to Reedy fibrancy).
 	\item Let $\cat_\infty$ $\infty$-category of small $\infty$-categories (\infref{item:examples}). In this case a Segal object is a {\it double $\infty$-category} in the sense of \cite[Definition 4.7]{haugseng2017highermorita} and \cite[Definition 4.1.7]{moser2020double} (up to Reedy fibrancy).
 	A {\it complete Segal object} in $\cat_\infty$ has not been studied independently, but is a first step towards the definition of a $2$-fold complete Segal space \cite{barwick2005nfoldsegalspaces}, which is a model of $(\infty,2)$-categories \cite{bergner2020surveyn,bergnerrezk2013comparisoni,bergnerrezk2020comparisonii}.
 	\item {\it Stable $\infty$-Category:} An $\infty$-category is stable if it has finite limits and colimits, the initial and final object coincide and the pushout and pullback square coincide \cite[Proposition 1.1.3.4]{lurie2017ha}, \cite[Theorem 4.4.12]{riehlverity2018elements}. We want to show that every complete Segal object $\kW$ in a stable $\infty$-category $\C$ is discrete. 
 	Indeed, combining \cref{prop:hoequivintomap} and \cref{def:cso} implies that $\kW_0 \to \kW_1$ is an inclusion, which means the square 
 	\begin{center}
 		\begin{tikzcd}
 			\kW_0 \arrow[r, "="] \arrow[d, "="']  \arrow[dr, phantom, "\ulcorner", very near start] & \kW_0 \arrow[d] \\
 			\kW_0 \arrow[r] & \kW_1 
 		\end{tikzcd}
 	\end{center}
 	is a pullback square and so, by stability, also a pushout square. As the map $\kW_0 \to \kW_1$ is now the pushout of the identity map $\id:\kW_0 \to \kW_0$, it is an equivalence and so the simplicial object $\kW$ is discrete.
 \end{enumerate}
 
\section{Two Applications: Strictification and Universality}
 In this last section we want to study very useful applications of our new found understanding of complete Segal objects and its interaction with Cartesian fibrations. The first one helps us construct complete Segal objects. The second one is about representing and hence better understanding the universal coCartesian fibration.
 
 \subsection{Building Segal Objects} \label{subsec:strictify}
 In this subsection we want to give a precise answer to the following problem (which was already outlined in \whyref{item:constructing cso}). 
 Assume we have a Segal Cartesian fibration $\R\to \C$ such that $\RFib(\R)_0$ and $\RFib(\R)_1$ are represented by $\kW_0,\kW_1$.
 How can we construct a Segal object $\widehat{\kW}$ representing $\R$?
 
 If we were in the setting of strict categories the answer would be straightforward: we take the fiber products and then use the Yoneda lemma several times. In the $\infty$-setting, however, we need a homotopy coherent approach and hence we will use strictification results for Cartesian fibrations.

 \begin{lemone} \label{the:building simp objects}
 	Let $\C$ be an $\infty$-category and  $\R \to \C$ be a Reedy right fibration that is level-wise representable. Then there exists a simplicial object $\kX: \DDelta^{op} \to \C$ such that 
 	we have an equivalence of Reedy right fibrations $\C_{/\kX} \simeq \R$.
 \end{lemone}
 
 \begin{proof}
  By \covref{item:left localization} $\R_\bullet$ is a simplicial $\infty$-category. 
  We can first construct the underlying complete Segal space $\Und(\R)_\bullet \to \Und(\C)$ (\infref{item:und}) and then apply $i_1^*$ \cite[Section 4]{joyaltierney2007qcatvssegal} to get a Reedy right fibration of quasi-categories $i_1^*\Und(\R)_\bullet \to \i_1^*\Und(\C)$ (analogous to \cite[Theorem 1.35]{rasekh2021cartfibmarkedvscso}).
  
  Applying the straightening construction \cite[Theorem 2.2.1.2]{lurie2009htt} level-wise (as discussed in \cite[Subsection 3.1]{rasekh2021cartfibmarkedvscso}), this fibration corresponds to a functor of simplicially enriched categories 
  $$F:\kC[i_1^*\Und(\C)]^{op} \to \Fun(\DDelta^{op},\sSet),$$
  where $\kC[i_1^*\Und(\C)]$ is the simplicially enriched category defined in \cite[Definition 1.1.5.1]{lurie2009htt}.
  Using the adjunction between products and functor categories this corresponds to a functor 
  $$\hat{F}_\bullet: \DDelta^{op} \to \Fun(\kC[i_1^*\Und(\C)]^{op},\sSet).$$
  The fact that the original fibration $\R \to \C$ is level-wise representable implies that each functor $\hat{F}_n$ is representable, meaning there is a natural equivalence of simplicially enriched functors $\hat{F}_n \simeq \Map(-,X_n)$. This means the functor $\hat{F}$ factors through the essential image of the Yoneda embedding
  $$\hat{F}_\bullet: \DDelta^{op} \to \kC[i_1^*\Und(\C)] \to \Fun(\kC[i_1^*\Und(\C)]^{op},\sSet).$$
  This gives us a simplicial object $\kX: \DDelta^{op} \to \C$, which gives us a representable Reedy right fibration $\C_{/\kX}$ (\cref{prop:rep left fib higher cat}/\ref{eq:contra table}).
  Finally, notice the straightening of $\C_{/\kX}$ is also equivalent to $\hat{F}$ and so $\C_{/\kX}$ and $\R$ are equivalent Reedy right fibrations. 
 \end{proof}
	
 \begin{lemone} \label{the:rep segal cart}
 	Let $\C$ be an $\infty$-category with finite limits and let $\R \to \C$ be a Segal Cartesian fibration.
 	Then $\R$ is representable if and only if $\RFib_0(\R)$ and $\RFib_1(\R)$ are representable.
 \end{lemone}
 
 \begin{proof}
 	One side is just a special case. So, let us assume $\RFib_0(\R)$ and $\RFib_1(\R)$ are representable right fibrations.
 	By \cref{the:building simp objects}, it suffices to prove that $\RFib_n(\R)$ is representable for $n \geq 2$, which means we have to show it has a final object. As $\R$ is a Segal Cartesian fibration (\reecovref{item:segal cart}) we have an equivalence of $\infty$-categories 
 	$$\RFib_n(\R) \to \RFib_1(\R) \underset{\RFib_0(\R)}{\times} ...\underset{\RFib_0(\R)}{\times} \RFib_1(\R).$$
 	Thus it suffices to prove that the right hand side has a final object. However, by the representability condition we have an equivalence of right fibrations
 	$$\RFib_1(\R) \underset{\RFib_0(\R)}{\times} ...\underset{\RFib_0(\R)}{\times} \RFib_1(\R) \simeq 
 	\C_{/\kW_1} \underset{\C_{/\kW_0}}{\times} ...\underset{\C_{/\kW_0}}{\times} \C_{/\kW_1}$$
 	where $\kW_1$ represents $\R_1$ and $\kW_0$ represents $\R_0$. Finally the right hand $\infty$-category has a final object if and only if the 
 	induced diagram of 
 	\begin{center}
 		\begin{tikzcd}
 			\kW_1 \arrow[r, "\kt"] & \kW_0 & \kW_1 \arrow[l, "\ks"'] \arrow[r] & ... & \kW_1 \arrow[l] \arrow[r, "\kt"]& \kW_0 & \kW_1 \arrow[l, "\ks"']
 		\end{tikzcd}
 	\end{center}
 	has a limit (by \covref{item:cocone}), which holds as $\C$ has finite limits. Here $\ks,\kt$ are the source target map $\ks,\kt: \kW_1 \to \kW_0$ (\infref{item:delta}). 
 \end{proof}
 
 \begin{theone} \label{the:segal Obj out of two obj}
 	Let $\kW_0$ and $\kW_1$ be two objects in an $\infty$-category $\C$ with finite limits.
 	Let $\R$ be a Segal Cartesian fibration over $\C$ such that $\R_0$ is represented by $\kW_0$ and $\R_1$ is represented by 
 	$\kW_1$. Then there exists a Segal object $\widehat{\kW}_{\bullet}$ such that $\widehat{\kW}_0 \simeq \kW_0$ and $\widehat{\kW}_1 \simeq \kW_1$.
 	Moreover, $\widehat{\kW}$ is complete if and only if $\R$ is a Cartesian fibration.
 \end{theone}
 
 \begin{proof}
 	By \cref{the:rep segal cart}, the Segal Cartesian fibration $\R$ is representable.
 	Thus, there exists a simplicial object $\widehat{\kW}$ such that $\R \simeq \C_{/\widehat{\kW}}$ and by universality we have $\widehat{\kW}_0 \simeq \kW_0$ and $\widehat{\kW}_1 \simeq \kW_1$. Moreover, $\R$ being a Segal Cartesian fibration implies that  $\widehat{\kW}$ is actually a Segal object (\cref{prop:segal obj vs Segal Cart}). Finally, again by \cref{prop:segal obj vs Segal Cart}, $\R$ is a Cartesian fibration if and only if $\widehat{\kW}$ is complete. 
 \end{proof}
  
 \subsection{Universal coCartesian Fibrations} \label{subsec:universal}
 In this final subsection we want to give another application: Understanding the {\it universal coCartesian fibration} via complete Segal objects (already mentioned in \whyref{item:universal fib}).
 This can also be found in \cite[Example 3.26]{stenzel2020comprehension}.
 Before we do so we want to review the analogous result for left fibrations. 
 
 Let $\C$ be an $\infty$-category and define $\newLFib_{/\C}$ as the $\infty$-category of left fibrations over $\C$ (given via the simplicial nerve \cite[Proposition 1.1.5.10]{lurie2009htt} applied to covariant model structure over $\C$ \infref{item:examples}). Then various authors (\cite[Theorem 2.2.1.2]{lurie2009htt}\cite{stevenson2017covariant}\cite{heutsmoerdijk2015leftfibrationi,heutsmoerdijk2016leftfibrationii}) have proven a natural equivalence of $\infty$-categories 
 \begin{equation}\label{eq:groth}
  \newLFib_{/\C} \simeq \s^\C.	
 \end{equation}
 Naturality implies that the {\it universal left fibration} is the left fibration over $\s$ that corresponds to the identity functor $\id: \s \to \s$. The problem is that the general construction of such fibrations is very abstract and difficult. We only know for certain that representable functors correspond to representable fibrations. Fortunately, in this case we can in fact represent the functor and compute the left fibration: The identity functor is equivalent to the functor represented by the point, $\map_\s(*,-): \s \to \s$. The left fibration constructed out of the representable functor is the under-category projection $\s_{*/} \to \s$, which more conventionally we think of as the $\infty$-category of pointed spaces $\s_* \to \s$. For a more thorough discussion of pointed spaces as the universal left fibration see \cite[Corollary 5.2.8]{cisinski2019highercategories}\cite[Lemma 2.2.4]{kazhdanvarshvsky2014yoneda} .
  
 We now want to generalize this result to {\it universal coCartesian fibrations}. Let $\coCartFib_{/\C}$ be the $\infty$-category of coCartesian fibration over $\C$ (\infref{item:examples})) and recall it is precisely the $\infty$-category of simplicial diagrams in $\newLFib_{/\C}$ that satisfy the complete Segal condition (\reecovref{item:cart}). Similarly let $\CSS$ be the $\infty$-category of complete Segal spaces and $\cat_\infty$, the $\infty$-category of (small) $\infty$-categories (again using the notational convention from \infref{item:examples}) and notice we still have an equivalence $\Und: \cat_\infty \to \CSS$ taking every $\infty$-category to its underlying $\infty$-category (\infref{item:und}). 
 
 If we take complete Segal objects in the $\infty$-category $\s^\C$ we get the $\infty$-category $\CSS^\C$. Indeed a simplicial diagram in space valued functors that satisfies the complete Segal conditions is precisely a complete Segal space valued functors. Applying these to the equivalence given in \ref{eq:groth} we get the equivalence 
 $$\coCartFib_{/\C} \simeq \CSS^\C \simeq (\cat_\infty)^\C$$
 where the second equivalence comes from $\Und^\C$. 
 (Notice an alternative, far more complicated, method for deducing the equivalence would have been to use the straightening construction for coCartesian fibrations \cite[Theorem 3.2.0.1]{lurie2009htt}).
 
 Again, the naturality implies that the {\it universal coCartesian fibration} is the coCartesian fibration over $\cat_\infty$ that corresponds to the identity functor $\id:\cat_\infty\to \cat_\infty$. Again, we cannot generally compute fibrations for arbitrary functors, however, we could compute the fibration if the identity map was representable. However, this is in fact false! There is no $\infty$-category $\C$ such that there is a natural equivalence $\map_{\cat_\infty}(\C,\D) \simeq \D$ for all $\infty$-categories $\D$ (as the left hand side is an $\infty$-groupoid). 
 The non-representability of this coCartesian fibration poses a serious challenge if we try to study it (as can be observed in \cite[Subsection 3.3.2]{lurie2009htt}).
 
 On the other hand, using the theory of representable coCartesian fibrations, we can in fact express the identity map as a representable functor. 
 Recall that $\cat_\infty$ has a cosimplicial object  $\DDelta \to \cat_\infty$ of the following form (\infref{item:fn}):
 \begin{center}
 	\begin{tikzcd}[row sep=0.5in, column sep=0.5in]
 		\1 
 		\arrow[r, shift left=1.2] \arrow[r, shift right=1.2]
 		& \2 
 		\arrow[l, shorten >=1ex,shorten <=1ex]
 		\arrow[r] \arrow[r, shift left=2] \arrow[r, shift right=2] 
 		& \3 
 		\arrow[l, shift right, shorten >=1ex,shorten <=1ex ] \arrow[l, shift left, shorten >=1ex,shorten <=1ex]
  		\arrow[r, shift right=1] \arrow[r, shift left=1] \arrow[r, shift right=3] \arrow[r, shift left=3]  		
 		& \cdots 
 		\arrow[l, shorten >=1ex,shorten <=1ex] \arrow[l, shift left=2, shorten >=1ex,shorten <=1ex] \arrow[l, shift right=2, shorten >=1ex,shorten <=1ex]	
 	\end{tikzcd}.
 \end{center}
 Moreover, recall the object $\mathbbe{n}$ represents the equivalence $\Und$ (\infref{item:fnund}), meaning for an $\infty$-category $\C$, we have an equivalence $\Map(\mathbbe{n},\C) \simeq \Und(\C)_n$. So $\Map(\mathbbe{n},-)$ represents the $n$-th space of the underlying complete Segal space. Hence, the functor represented by this cosimplicial object $\Diag \circ\Map(\mathbbe{n},-): \cat_\infty \to \cat_\infty$ is equivalent to the identity map.
 This proves that the {\it universal coCartesian fibration} is given by the representable coCartesian fibration $(\cat_\infty)_{\mathbbe{n}/}\to \cat_\infty$.
  
\bibliographystyle{alpha}
\bibliography{main}

\end{document}